\documentclass{svmult}

\usepackage{makeidx}         
\usepackage{graphicx}        
\usepackage{multicol}        
\usepackage[bottom]{footmisc}
\makeindex             

\newcommand{\EssVar}{\mbox{EssVar}}
\newcommand{\NEssVar}{N_{\mbox{ess}}}

\begin{document}

\title*{Reducing the number of variables of a polynomial}
\author{Enrico Carlini}
\institute{Politecnico di Torino, Corso Duca degli Abruzzi 24,
10129 Torino, Italia \texttt{enrico.carlini@polito.it}}
%
%
\maketitle

\section{Abstract}
\label{sec:1}

In this paper, we consider two basic questions about presenting a
homogeneous polynomial $f$: how many variables are needed for
presenting $f$? How can one find a presentation of $f$ involving
as few variables as possible? We give a complete answer to both
questions, determining the minimal number of variables needed,
$\NEssVar(f)$, and describing these variables through their linear
span, $\EssVar(f)$. Our results give rise to effective algorithms
which we implemented in the computer algebra system {\bf CoCoA}
\cite{cocoa}.

\section{Introduction}
\label{sec:2}

Polynomials, also seen as symmetric tensors, are ubiquitous in
Applied Mathematics. They appear in Mechanics (\cite{McL}), Signal
and Image Processing (\cite{ComMour}), Algebraic Complexity Theory
(\cite{BuClSh}), Coding and Information Theory (\cite{Ro}), etc..

One of the main open issue is to manipulate polynomials in order
to obtain presentations suiting the special needs of the
application at hand.

In Mechanics, it is often useful to {\em separate variables}.
Given a polynomial $f(x_1,\ldots,x_n)$, one splits the set of
variables in two pieces, e.g. $\{x_1,\ldots,x_r\}$ and
$\{x_{r+1},\ldots,x_n\}$, and a presentation of $f$ is searched of
the following type
\[
f(x_1,\ldots,x_n)=g(x_1,\ldots,x_r)+h(x_{r+1},\ldots,x_n)
\]
for some polynomials $g$ and $h$.

Separating variables is a well established technique and the
search for splitting methods in general is very active (see
\cite{McL}).

In Signal Processing, homogeneous polynomials (also known as {\em
quantics} from ancient Invariant Theory) are of crucial
importance. The main interest is in the so called {\em sum of
powers presentations}, where a homogeneous polynomial $f$ of
degree $d$ is presented as
\[
f=l_1^d+\ldots +l_s^d
\]
where $l_1,\ldots,l_s$ are linear forms.

Sum of powers presentations are treated in connection with
quantics in \cite{ComMour}, while a more general approach relating
them to Polynomial Interpolation and Waring Problem can be found
in \cite{Ci01}.

In this paper, we consider two basic questions about presenting a
homogeneous polynomial (from now on referred to as a {\em form})
in a ``easier'' way. Given a form $f$, how many variables are
needed for presenting it? How can one find a presentation of $f$
involving as few variables as possible?

Even if these problems are so natural, we are not aware of a
complete solution existing in the literature. In this paper, we
give a complete answer to both questions. Our results give rise to
effective algorithms which we implemented in the computer algebra
system {\bf CoCoA} (freely available at {\tt
cocoa.dima.unige.it}).

More precisely, given a form $f\in S=k[x_1,\ldots,x_n]$, $k$ any
field, we call {\em essential number of variables of $f$} the
smallest integer $r$ for which there exists a set of linear forms
$\{y_1,\ldots,y_r\}\subset S$ such that
\[
f\in k[y_1,\ldots,y_r];
\]
the linear forms $y_1,\ldots,y_r$ are called {\em essential
variables of $f$}. Then our main result is (see Definition
\ref{catDEF}, Definition  \ref{essentialDEF} and Section
\ref{apolaritySEC} for the notation involved):

\medskip\noindent
{\bf Proposition \ref{propNEssVar&EssVar}\ }{\em Let $f$ be a
homogeneous element in $S=k[x_1,\ldots,x_n]$ and
$T=k[\partial_1,\ldots,\partial_n]$ denote the ring of
differential operators. Then
\[
\NEssVar(f)=\mbox{rk}(\mathcal{C}_f),
\]
i.e. the number of essential variables of $f$ is the rank of its
first catalecticant matrix, and
\[
EssVar(f)=\langle D\circ f : D \in T_{d-1}\rangle,
\]
i.e. the essential variables of $f$ span the space of its
$(d-1)^{\mbox{th}}$ partial derivatives.}
\medskip

In Section \ref{apolaritySEC}, we briefly recall some facts from
Apolarity Theory which are the main tools of our analysis.

In Section \ref{howmanySEC}, we use Apolarity and Catalecticant
Matrices to obtain our main result. In Subsection \ref{pcSEC}, we
give some examples of the use of our algorithms.

\begin{remark}
In this paper we work with forms, i.e. homogeneous polynomials. To
apply our results to {\em any} polynomial $f$, it is enough to
work with its homogenization $f^h$ with respect to a new variable.
Clearly, e.g., a presentation of $f^h$ in essential variables
readily produces a presentation of $f$ in essential variables: it
is enough to dehomogenize.
\end{remark}

\begin{remark}\label{fieldRM}
Throughout the paper $k$ will denote a field of characteristic 0.
Our results also hold in positive characteristic, but more
advanced techniques are required for proving them, e.g.
differentiations have to be replaced with contractions and divided
powers have to be introduced (see \cite{IaKa}).
\end{remark}

The author wishes to thank B. Reznick and C. Ciliberto for their
ideas on the problem. The {\bf CoCoA} Team in Genoa, and
especially Anna Bigatti, were of great help in the implementation
of the algorithms. The comments and criticisms of the two
anonymous referees were of help in improving the presentation of
the results.

\section{Apolarity}\label{apolaritySEC}

In this section we will briefly recall some basic facts from {\em
Apolarity Theory} or, in modern terms, {\em Inverse Systems
Theory}. Comprehensive references are \cite{Ge}, mainly Lecture 2,
6 and 8, and \cite{IaKa}.

Consider the polynomial rings
\[
S=k[x_1,\ldots,x_n]\mbox{ and }T=k[\partial_1,\ldots,\partial_n],
\]
where $k$ denotes a field of characteristic 0, and give $S$ a
$T$-module structure via differentiation, i.e. we will think of
$T$ as the ring of differential operators acting on $S$. We denote
this action with ``$\circ$'', e.g. $\partial_j\circ
f=\frac{\partial}{\partial x_j} f$ for $f\in S$.

There is a natural perfect paring between homogeneous pieces of
the same degree of $S$ and $T$, namely
\[
\begin{array}{ccccl}
S_i & \times & T_i & \longrightarrow & k \\
f   &        & D   &                 & D\circ f
\end{array},
\]
is a perfect pairing for all $i$; in particular, $S_i$ and $T_i$
are dual to each other. Given subspaces,
\[
V\subseteq S_i\mbox{ and } W\subseteq T_i
\]
we denote by
\[
V^\perp\subseteq T_i\mbox{ and } W^\perp\subseteq S_i
\]
their orthogonal with respect to this pairing; notice, e.g., that
$\dim_k V + \dim_k V^\perp =\dim_k S_i=\dim_k T_i$.

Given a form $f\in S_d$, the ideal
\[
f^\perp=\{D\in T : D\circ f=0\}
\]
is a homogeneous ideal of $T$ and it is called the {\em
orthogonal} ideal of $f$.


Orthogonal ideals play a central role in the theory: they contain
all the differential operators annihilating a given form and even
more information, as it is shown by the following Lemma (for a
proof see \cite{Ge}, Proposition 8.10).
\begin{lemma}\label{saturatedLEM} Let $f$ be a degree $d$ form in $S$, then $D\in
T_i, i<d$, is such that
\[
D\circ f=0
\]
if and only if
\[
D\circ(D'\circ f)=0
\]
for all $D'\in T_{d-i}$. In other terms, for $0<i<d$,
$(f^\perp)_i$ is orthogonal to the $k$-vector space spanned by the
$(d-i)^{\mbox{th}}$ partial derivatives of $f$.
\end{lemma}

Orthogonal ideals can be easily described introducing {\em ad hoc}
matrices. In this paper it will be enough to describe the degree
one part of a given orthogonal ideal, but similar descriptions
exist in each degree.

\begin{definition}\label{catDEF}
Let $f\in S_d$ and fix the standard monomial basis, e.g. with
respect to lex order, $\{M_1,\ldots, M_N\}$ of the $k$-vector
space $S_{d-1}$. For $i=1,\ldots,n$, consider the first partials
\[
\partial_i\circ f=c_{i1}M_1+\ldots +c_{iN} M_N.
\]
The {\em first catalecticant matrix of $f$} is
\[
\left(\mathcal{C}_f\right)_{ij}=c_{ij},
\]
$i=1,\ldots,n, j=1,\ldots,N$.
\end{definition}
For a general treatment of catalecticant matrices and their
applications see \cite{ge1}, \cite{Ge} and \cite{IaKa}.
\begin{example}\label{ex1}
Let $f=x_1x_2x_3\in k[x_1,x_2,x_3]$ and consider the monomial
basis
\[
\lbrace \ x_{1\, }^{2}, x_{1\, }x_{2\, }, x_{1\, }x_{3\, }, x_{2\,
}^{2}, x_{2\, }x_{3\, }, x_{3\, }^{2}\ \rbrace
\]
of the space of degree two forms. Then
\[
\mathcal{C}_f=\left( \begin{array}{llllll}
0 & 0 & 0 & 0 & 1 & 0 \\
0 & 0 & 1 & 0 & 0 & 0 \\
0 & 1 & 0 & 0 & 0 & 0 \end{array}\right).
\]
\end{example}

Catalecticant matrices determine the degree one part of orthogonal
ideals readily:
\begin{lemma}\label{catalectLEM}
Let $f\in k[x_1,\ldots,x_n]$ be a form, then
\[
(a_1\partial_1+\ldots +a_n\partial_n)\circ f=0
\]
if and only if the vector $(a_1,\ldots,a_n)$ is in the left kernel
of $\mathcal{C}_f$. In particular,
$\dim_k(f^\perp)_1=n-\mbox{rk}(\mathcal{C}_f)$.
\end{lemma}
\begin{proof}
The statement simply follows writing down the action of
$a_1\partial_1+\ldots +a_n\partial_n$ on $f$ componentwise and
considering the corresponding linear system of equations.
\end{proof}
\begin{remark}\label{rk1rem}
Let $l\in k[x_1,\ldots,x_n]$ be a linear form and consider its
$d$-th power $f=l^d$. Let $L$ be a linear differential operator
and notice that $L\circ f=0$ if and only if $L\circ l=0$ which is
a linear equation in the coefficients of $L$. Hence
$\mbox{rk}(\mathcal{C}_{l^d})=1$ (actually, even the converse is
true). In particular, this means that the form of Example
\ref{ex1} is not a pure power.
\end{remark}

\section{How many variables?}\label{howmanySEC}

In this section we will use apolarity to answer our two basic
questions: how many variables do we need to present a given form?
How can we find a presentation involving as few variables as
possible?

In what follows, we will work with the polynomial ring
$S=k[x_1,\ldots,x_n]$, where $k$ is {\em any} field such that
$\mbox{char}(k)=0$ (in positive characteristic similar results
hold, but, in this paper, we decided to avoid the technical
difficulties involved).

Lets introduce some definitions:
\begin{definition}\label{essentialDEF}
Given a form $f$ in $S$, the {\em number of essential variables of
$f$}, $\NEssVar(f)$, is the smallest integer $r$ such that there
exist linear forms $y_1,\ldots,y_r\in S$ for which $f\in
k[y_1,\ldots, y_r]$. We call {\em essential variables of $f$} any
set of generators of the $k$-vector space $\EssVar(f)=\langle
y_1,\ldots, y_r\rangle$.
\end{definition}
Roughly speaking, given a form $f\in S$, $\NEssVar(f)$ tells us
how many variables are necessary for presenting $f$, while
$\EssVar(f)$ tells us how we can find such variables. In
particular, it is clear that, if
\[
\NEssVar(f)=r\mbox{ and }\EssVar(f)=\langle y_1,\ldots,y_r
\rangle,
\]
then there exists $g\in k[y_1,\ldots, y_r]\subset S$ such that
$f=g$.
\begin{example}
Consider the form $f=f(x_1,x_2,x_3)=(x_1+x_2){(x_1-x_3)}^2$ in
$k[x_1,x_2,x_3]$. Clearly $f$ is an element of the subring
$k[y_1,y_2]$, where $y_1=x_1+x_2$ and $y_2=x_1-x_3$. Hence
$\NEssVar(f)\leq 2$ and equality holds by Remark \ref{rk1rem}, as
$\mbox{rk}(\mathcal{C}_f)\neq 1$ and $f$ is not a pure power.
Also, notice that $\EssVar(f)=\langle x_1+x_2, x_1-x_3\rangle$ and
sets of possible essential variables are: $\{x_1+x_2,x_1-x_3\}$,
$\{x_2+x_3,2x_1+x_2-x_3\}$, etc..
\end{example}

Using apolarity we can effectively determine $\NEssVar$ and
$\EssVar$ for a given form:
\begin{proposition}\label{propNEssVar&EssVar}
Let $f$ be a homogeneous element in $S=k[x_1,\ldots,x_n]$ and
$T=k[\partial_1,\ldots,\partial_n]$ denote the ring of
differential operators. Then
\[
\NEssVar(f)=\mbox{rk}(\mathcal{C}_f),
\]
i.e. the number of essential variables of $f$ is the rank of its
first catalecticant matrix, and
\[
EssVar(f)=\langle D\circ f : D \in T_{d-1}\rangle,
\]
i.e. the essential variables of $f$ span the space of its
$(d-1)^{\mbox{th}}$ partial derivatives.
\end{proposition}
\begin{proof}
If $\NEssVar(f)=r$, then $f\in k[y_1,\ldots,y_r]$ for some linear
forms $y_1,\ldots, y_r$ in $S$. Let
\[
\langle y_1,\ldots, y_r\rangle^\perp=\langle L_1,\ldots,
L_{n-r}\rangle\subset T_1
\]
and notice that $(f^\perp)_1\supseteq \langle L_1,\ldots,
L_{n-r}\rangle$. Thus, by Lemma \ref{catalectLEM}, we have
$\mbox{rk}(\mathcal{C}_f)\leq\NEssVar(f)$.


If $\mbox{rk}(\mathcal{C}_f)=t$, then $(f^\perp)_1=\langle
D_1,\ldots, D_{n-t}\rangle$. Complete this to a basis of $T_1$
\[
\langle D_1,\ldots, D_{n-t},Y_1,\ldots,Y_t\rangle
\]
and consider the dual basis of $S_1$ defined by the apolarity
perfect pairing
\[
\langle z_1,\ldots, z_{n-t},y_1,\ldots,y_t\rangle.
\]
Hence, after a linear change of variables, we have
$f=f(z_1,\ldots, z_{n-t},y_1,\ldots,y_t)$. But $D_j$ annihilates
all the elements of the chosen basis of $S_1$ but $z_j$. As
$(f^\perp)_1=\langle D_1,\ldots, D_{n-t}\rangle$ we conclude that
\[
f\in k[y_1,\ldots,y_t]
\]
and $\mbox{rk}(\mathcal{C}_f)\geq\NEssVar(f)$.

To conclude the proof, notice that the prefect pairing $S_1\times
T_1\rightarrow k$ induces a well defined perfect pairing of
$k$-vector spaces
\[
V\times\left(\frac{T}{f^\perp}\right)_1\longrightarrow k
\]
where
\[
V=\left((f^\perp)_1\right)^\perp=\langle l : l\in S_1, L\circ
l=0\mbox{ for all } L\in(f^\perp)_1\rangle
\]
and, with the notations above,
$\left(\frac{T}{f^\perp}\right)_1=\langle Y_1,\ldots,Y_t\rangle$
and hence $V=\EssVar(f)$. The result follows applying Lemma
\ref{saturatedLEM} ($i=1$ case) which yields
\[
V=\langle D'\circ f: D'\in T_{d-1} \rangle.
\]
\end{proof}

\begin{example}\label{ex2}
Given the form
\[
f=x_{1\, }^{3} + x_{1\, }^{2}x_{2\, }-2x_{1\, }^{2}x_{3\,
}-2x_{1\, }x_{2\, }x_{3\, } + x_{1\, }x_{3\, }^{2} + x_{2\,
}x_{3\, }^{2}
\] we want to determine $\NEssVar(f)$ and
$\EssVar(f)$. In order to apply Proposition
\ref{propNEssVar&EssVar}, we compute the first catalecticant
matrix of $f$
\[\mathcal{C}_f=
\left( \begin{array}{rrrrrr}
3 & 2 & -4 & 0 & -2 & 1 \\
1 & 0 & -2 & 0 & 0 & 1 \\
-2 & -2 & 2 & 0 & 2 & 0 \end{array}\right).
\]
Hence $\NEssVar(f)=\mbox{rk}(\mathcal{C}_f)=2$ and $f$ can be
presented as a form in two variables. To determine the essential
variables of $f$, it is enough to compute the span of the second
partial derivatives of $f$:
\[\EssVar(f)=\langle x_2+x_3,x_1-x_3\rangle.\]
Summing these up, we see that there exists a degree 3 form
$g(y_1,y_2)\in k[y_1,y_2]$ such that
\[g(x_2+x_3,x_1-x_3)=f(x_1,x_2,x_3),\]
but how can we find $g$?
\end{example}

To complete our analysis, we want to present a form $f$ as a
polynomial only involving essential variables: this can be done
almost tautologically, but the notation are quite involved. We
begin with an example.
\begin{example}\label{ex3}
Consider the form $f\in S=k[x_1,x_2,x_3]$ in Example \ref{ex2}. We
already showed that there exists $g\in k[y_1,y_2]\subset S$ such
that $f=g$. To determine $g(y_1,y_2)$, consider
$\EssVar(f)=\langle x_2+x_3,x_1-x_3\rangle$ and complete its basis
to a basis of $S_1$: we choose
$\{y_1=x_2+x_3,y_2=x_1-x_3,z_1=x_1\}$. Hence we have a linear
change of variables given by
\[
\left\lbrace
\begin{array}{l}
x_1=z_1,\\
x_2=y_1+y_2-z_1,\\
x_3=z_1-y_2.
\end{array}
\right.
\]
The basic requirement of the form $g(y_1,y_2)$ is to satisfy the
relation
\[
g(x_2+x_3,x_1-x_3)=f(x_1,x_2,x_3).
\]
From this, changing variables, we get
\[
g(y_1,y_2)=f(z_1,y_1+y_2-z_1,z_1-y_2)=y_1y_2^2+y_2^3,
\]
which is the desired presentation in essential variables. As a
byproduct, we readily see that
\[
f=(x_2+x_3)(x_1-x_3)^2+(x_1-x_3)^3
\]
which is quite surprising considering the original presentation
\[
f=x_{1\, }^{3} + x_{1\, }^{2}x_{2\, }-2x_{1\, }^{2}x_{3\,
}-2x_{1\, }x_{2\, }x_{3\, } + x_{1\, }x_{3\, }^{2} + x_{2\,
}x_{3\, }^{2}.
\]
\end{example}

The procedure showed in the previous Example works in general.
Given a form $f=f(x_1,\ldots,x_n)\in S$, we compute
$\NEssVar(f)=r$ and we choose a basis for $\NEssVar(f)=\langle
y_1,\ldots, y_r\rangle$; to avoid triviality, assume $r<n$. Now,
our goal is to determine $g=g(y_1,\ldots,y_r)\in
k[y_1,\ldots,y_r]\subset S$ such that $f=g$. To do this, complete
the basis of $\NEssVar(f)\subset S_1$ to a basis of $S_1$
\[
S_1=\langle y_1,\ldots, y_r, z_1,\ldots,z_{n-r}\rangle.
\]
As $S_1=\langle x_1,\ldots,x_n \rangle$, the completed basis
yields a linear change of variables
\[(\dagger)\left\lbrace
\begin{array}{l}
x_1=x_1(y_1,\ldots, y_r, z_1,\ldots,z_{n-r}),\\
\vdots\\
x_n=x_n(y_1,\ldots, y_r, z_1,\ldots,z_{n-r}).
\end{array}
\right.
\]
Notice that $y_1,\ldots,y_r$ are linear forms in $S$ and hence
there exist linear functions such that
$y_i=y_i(x_1,\ldots,x_r),i=1,\dots,r$. Moreover, the following
identities hold by construction of $(\dagger)$
\[
y_i=y_i(x_1(y_1,\ldots, y_r,
z_1,\ldots,z_{n-r}),\ldots,x_r(y_1,\ldots, y_r,
z_1,\ldots,z_{n-r}))
\]
for $i=1,\ldots ,n$.

To determine $g$, it is enough to consider the desired relation
\[
f(x_1,\ldots,x_n)=g(y_1(x_1,\ldots,x_r),\ldots,y_r(x_1,\ldots,x_r)).
\]
and to apply the linear change of variables $(\dagger)$. Thus we
obtain $g(y_1,\ldots,y_r)$:


\[\begin{array}{c}
g(y_1,\ldots,y_r) = \\

=g(y_1(x_1(y_1,\ldots, y_r,
z_1,\ldots,z_{n-r}),\ldots,x_r(y_1,\ldots, y_r,
z_1,\ldots,z_{n-r})),\ldots \\
\ldots,y_r(x_1(y_1,\ldots, y_r,
z_1,\ldots,z_{n-r}),\ldots,x_r(y_1,\ldots, y_r,
z_1,\ldots,z_{n-r})))= \\
=f(x_1(y_1,\ldots, y_r,
                  z_1,\ldots,z_{n-r}),\ldots,x_n(y_1,\ldots, y_r,
                  z_1,\ldots,z_{n-r})).
\end{array}
\]
Notice that, as $f$ and the functions $x_i(y_1,\ldots, y_r,
z_1,\ldots,z_{n-r}),i=1,\ldots,n,$ are {\em explicitly} known, we
have completely determined $g$ as an element in
$k[y_1,\ldots,y_r]$.
\begin{remark}\label{cylinderRM}
As a straightforward application of the theory, we consider the
detection of cylinders (i.e. algebraic surfaces ruled by a family
of parallel lines moving along a fixed curve). Suppose you are
given the polynomial equation of a surface
$\mathcal{F}:f(x,y,z)=0$ in three space and you want to decide
whether $\mathcal{F}$ is a cylinder or not. It is well known that
$\mathcal{F}$ is a cylinder if and only if its defining equation
is a function of two planes, i.e. there exist linear forms
$m(x,y,z)$ and $l(x,y,z)$ such that $f(x,y,z)=g(m,n)$ for some
polynomial $g$. Hence, we readily have an effective procedure for
cylinder detection:
\[
\mathcal{F} \mbox{ is a cylinder if and only if }\NEssVar(f^h)\leq
3,
\]
where $f^h$ denotes the homogenization of $f$ (see Example
\ref{cylinderEX}). Clearly, the method applies in any dimension
for deciding whether a given hypersurface is a cylinder or not.
\end{remark}

\subsection{Using a computer}\label{pcSEC}

The results of our analysis can be easily translated into
algorithms and we wrote down procedures to be used with the
Computer Algebra system {\bf CoCoA}.

We begin with reporting a {\bf CoCoA} session illustrating the use
of our algorithms to work out the expository Examples \ref{ex2}
and \ref{ex3}.
\begin{example}
First we define the form we want to study
\begin{verbatim}
F:=x^3 + x^2y - 2x^2z - 2xyz + xz^2 + yz^2;
\end{verbatim}

To compute the number of essential of variables, use the function
{\tt NEssVar(F)}:
\begin{verbatim}
NEssVar(F);
2
-------------------------------
\end{verbatim}
To determine a choice of essential variables, use the function
{\tt EssVar(F)}:
\begin{verbatim}
EssVar(F);
[y + z, x - z]
-------------------------------
\end{verbatim}
Finally, {\tt NewPres(F)} produces a presentation of the form
involving the essential variables {$\mathtt{y[1]=y+z,y[2]=x-z}$}:
\begin{verbatim}
NewPres(F);
y[1]y[2]^2 + y[2]^3
-------------------------------
\end{verbatim}
\end{example}

Usually, a given polynomial $f(x_1,\ldots,x_n)$ will essentially
involve $n$ variables, i.e. $\NEssVar(f)=n$. Hence our algorithms
{\em do not} help in solving the polynomial equation $f=0$.
Nevertheless, our procedure should be used as a pre-processing
tool. In fact, {\em if} the number of variables can be decreased,
then the numerical solution of the equation can be performed much
more efficiently. We illustrate this with the following
``extreme'' example.
\begin{example}
We consider the degree three polynomial in four variables
\[
f(x,y,z,t)=f_0(x,y,z,t)+f_1(x,y,z,t)+f_2(x,y,z,t)+f_3(x,y,z,t),
\]
where
\[
\begin{array}{lcl}
f_0 & = & 3 \\
f_1 & = & - x - y + 2z + 3t\\
f_2 & = & 5 x^{2} + 10xy + 5y^{2}-20xz-20yz + 20z^{2}-30xt-30yt +
60zt + 45t^{2} \\
f_3 & = & x^{3} + 3x^{2}y + 3xy^{2} + y^{3}-6x^{2}z-12xyz-6y^{2}z
+ 12xz^{2} + 12yz^{2}+\\
    &  & -8z^{3}-9x^{2}t-18xyt-9y^{2}t + 36xzt + 36yzt-36z^{2}t + 27xt^{2} +
    27yt^{2}+ \\
    &  & -54zt^{2}-27t^{3}.
\end{array}
\]
In order to solve the equation $f(x,y,z,t)=0$, we apply our
algorithms to the degree 2 and 3 pieces of $f$:
\begin{verbatim}
EssVar(F2);
[x + y - 2z - 3t]
-------------------------------
NewPres(F2);
5y[1]^2
-------------------------------
\end{verbatim}
and hence $f_2(x,y,z,t)=5{y_1}^2$, where $y_1=x+y-2z-3t$.
Similarly
\begin{verbatim}
EssVar(F3);
[x + y - 2z - 3t]
-------------------------------
NewPres(F3);
y[1]^3
-------------------------------
\end{verbatim}
which yields $f_3(x,y,z,t)={y_1}^3$. In conclusion, to solve the
equation $f(x,y,z,t)=0$, it is enough to solve the equation in one
variable
\[
{y_1}^3+5{y_1}^2-y_1+3=0
\]
and to apply some linear algebra to find all the solutions.
\end{example}

We conclude with a Geometric example about cylinder detection.
\begin{example}\label{cylinderEX}
Consider the degree five surface in three space $\mathcal{F}:
f(x,y,z)=0$, where
\[
f=f_0+f_2+f_5
\]
and
\[
f_0=-1, f_2=x^{2} - xy-2y^{2}-3yz - z^{2},
\]
\[
\begin{array}{lll}
f_5 & = &x^{5} +
2x^{4}y-2x^{3}y^{2}-8x^{2}y^{3}-7xy^{4}-2y^{5}+3x^{4}z-18x^{2}y^{2}z-24xy^{3}z+\\
    &   &-9y^{4}z + 2x^{3}z^{2}-12x^{2}yz^{2}-30xy^{2}z^{2}-16y^{3}z^{2}-2x^{2}z^{3}-16xyz^{3}+\\
    &   &-14y^{2}z^{3}-3xz^{4}-6yz^{4}- z^{5}.
\end{array}
\]
In order to decide whether $\mathcal{F}$ is a cylinder or not, we
follow Remark \ref{cylinderRM}. Introduce a new variable $t$ and
consider the homogenization of $f$, $f^h=t^5 f_0+t^3 f_2+f_5$.
Using {\bf CoCoA} and denoting by {\tt FH} the form
$f^h(x,y,z,t)$, we get:
\begin{verbatim}
NEssVar(FH);
3
-------------------------------
EssVar(FH);
[t, y + 2/3z, x + 1/3z]
-------------------------------
NewPres(FH);
-y[1]^5 - 2y[1]^3y[2]^2 - 2y[2]^5 - y[1]^3y[2]y[3] -
7y[2]^4y[3] + y[1]^3y[3]^2 - 8y[2]^3y[3]^2 - 2y[2]^2y[3]^3 +
2y[2]y[3]^4 + y[3]^5
-------------------------------
\end{verbatim}
In conclusion, $f^h(x,y,z,t)=g(y_1,y_2,y_3)$ where $g$ is the
output of the function {\tt NewPres(FH)} and
\[
\left\lbrace
\begin{array}{lll}
y_1=t\\
\\
y_2=y+\frac{2}{3}z\\
\\
y_3=x+\frac{1}{3}z
\end{array}
\right. .\]

Hence, we have the polynomial equality $f(x,y,z)=g(1,y_2,y_3)$ and
$\mathcal{F}$ is a cylinder ruled by lines parallel to the line
$y_2=y_3=0$.
\end{example}

\bibliographystyle{alpha}
\bibliography{carlini}



\printindex
\end{document}